\documentclass[12pt]{amsart}
\usepackage{amsfonts}
\usepackage{ifthen}
\usepackage{amsthm}
\usepackage{amsmath}
\usepackage{graphicx}
\usepackage{amscd,amssymb,amsthm}
\usepackage{graphicx}
\usepackage{epstopdf}
\usepackage{hyperref}
\usepackage{color}

\newcounter{minutes}
\divide\time by 60
\newcounter{hours}
\multiply\time by 60 \addtocounter{minutes}{-\time}

\newtheorem{lemma}{Lemma}[section]
\newtheorem{theorem}{Theorem}[section]
\newtheorem{corollary}{Corollary}[section]
\newtheorem{remark}{Remark}[section]

\newtheorem{example}{Example}[section]

\newcommand{\real}{\operatorname{Re}}


\title{Partial sums of Hyper-Bessel function with applications}

\author[\.{I}. Akta\c{s}]{\.{I}brah\.{I}m Akta\c{s}}
\address{Department of Mathematical Engineering, Faculty of Engineering and Natural Sciences, G\"{u}m\"{u}\c{s}hane University, G\"{u}m\"{u}\c{s}hane, Turkey}
\email{aktasibrahim38@gmail.com}

\begin{document}

\maketitle

\begin{abstract}
The main purpose of present paper is to determine some lower bounds for the quotient of the normalized hyper-Bessel function and its partial sum, as well as for the quotient of the derivative of normalized hyper-Bessel function and its partial sum. In addition, some applications related to obtained results are also given.
\end{abstract}

\section{Introduction and Preliminaries}
There are a vivid interest on special and geometric function theories  since their some close relations. Actually, there are various developments regarding partial sums of analytic univalent functions in the recent years. The readers may found these interesting developments in the papers \cite{AO1,Frasin,OSS,Sheil,Silverman,Silvia}. Also, partial sums of some special functions and their applications have been considered by the authors in \cite{AO2,AO3,CD,OY,RADU,YO}. Especially, the authors in \cite{OY} investigated partial sums of generalized Bessel functions first time for the special functions in 2014.  As a result of this, this problem has been carried to other special functions. Later, many authors have investigated the same problem for other special functions like Struve, Lommel, Mittag-Leffler, $ q $-Bessel and Dini functions. Motivated by the previous works on analytic univalent functions and special functions our main aim is to determine some lower bounds for the quotient of normalized hyper-Bessel functions and its partial sum, as well as for the quotient of the derivative of normalized hyper-Bessel function and its partial sum. In addition, we give some applications regarding our main results.

Before starting our main results we would like give some basic concepts concerning geometric function theory and the definition of hyper-Bessel function which is a natural extentions of classical Bessel function of the first kind.

Let $\mathcal{A}$ denote the class of functions of the following form: 
\begin{equation}\label{eq1.1}
f(z)=z+\sum_{n\geq2}a_{n}z^{n},
\end{equation}
which are analytic in the open unit disk  $$\mathcal{U}=\{z:z\in\mathbb{C}\text{ and }\left|z\right|<1\}.$$
We denote by $\mathcal{S}$ the class of all functions in  $\mathcal{A}$ which are univalent in $\mathcal{U}$.

Now, consider the hyper-Bessel function defined as follow: (see \cite{Chaggara})
\begin{equation}\label{eq1.2}
J_{\alpha_d}(z)=\dfrac{\left(\frac{z}{d+1}\right)^{\alpha_1+\dots+\alpha_d}}{\prod_{i=1}^{d}\Gamma\left(\alpha_i+1\right)}{_0F_d}
\left(\begin{matrix}-\\(\alpha_{d}+1)\end{matrix};-\left(\frac{z}{d+1}\right)^{d+1}\right),
\end{equation}
where the notation
\begin{equation}\label{eq1.3}
{_pF_q}\left(\begin{matrix}(\beta_p)\\(\gamma_q)
\end{matrix};x\right)=\sum_{n\geq0}\frac{(\beta_1)_n(\beta_2)_n\cdots(\beta_p)_n}{(\gamma_1)_n(\gamma_2)_n\cdots(\gamma_p)_n}\frac{x^n}{n!}
\end{equation}
denotes the generalized hypergeometric function,
$(\beta)_n$ is the shifted factorial
(or Pochhammer's symbol) defined by $(\beta)_0=1, (\beta)_n=\beta(\beta+1)\cdots(\beta+n-1), n\geq1$ and the contracted notation $\alpha_d$ is used to abbreviate the array of $d$ parameters $\alpha_1,\dots, \alpha_d$.

By considering the equations \eqref{eq1.2} and \eqref{eq1.3} it can be easily seen that the function $z\mapsto J_{\alpha_d}(z)$ has the following infinite sum representation
\begin{equation}\label{eq 1.4}
J_{\alpha_d}(z)=\sum_{n\geq0}\frac{(-1)^n}{n!\prod_{i=1}^{d}\Gamma\left(\alpha_i+1+n\right)}\left(\frac{z}{d+1}\right)^{n(d+1)+\alpha_1+\dots+\alpha_d}.
\end{equation}

The normalized hyper-Bessel function $\mathcal{J}_{\alpha_d}(z)$ is defined by
\begin{equation}\label{eq 1.5}
J_{\alpha_d}(z)=\frac{\left(\frac{z}{d+1}\right)^{\alpha_1+\dots+\alpha_d}}{\prod_{i=1}^{d}\Gamma\left(\alpha_i+1\right)}\mathcal{J}_{\alpha_d}(z).
\end{equation}
By combining the equations \eqref{eq 1.4} and \eqref{eq 1.5} we get the following series representation
\begin{equation}\label{eq1.6}
\mathcal{J}_{\alpha_d}(z)=\sum_{n\geq0}\frac{(-1)^n}{n!\prod_{i=1}^{d}(\alpha_i+1)_n}\left(\frac{z}{d+1}\right)^{n(d+1)}.
\end{equation}
Since the function $ \mathcal{J}_{\alpha_d} $ does not belong to the class $\mathcal{A}$, we consider the following normalized form
\begin{equation}\label{eq1.7}
{f_{\alpha_d}(z)}=z\mathcal{J}_{\alpha_d}(z)=\sum_{n\geq0}A_nz^{n(d+1)+1},
\end{equation}
where $A_n=\frac{(-1)^n}{n!\left(d+1\right)^{n(d+1)}\prod_{i=1}^{d}(\alpha_i+1)_n}.$
As consequence of this consideration we have that the function ${f_{\alpha_d}}\in\mathcal{A}$. 
Here, we would like to mention that the following inequalities
\begin{equation}\label{eq1.8}
n!\geq2^{n-1}
\end{equation}
and
\begin{equation}\label{eq1.9}
\left(\alpha_i+1\right)^n\geq\left(\alpha_i+1\right)_n
\end{equation}
hold true for $n\in\mathbb{N}=\{1,2,\dots\}$ and $i\in\{1,2,\dots,d\},$ which will be used in the proof of our main results. Also, we will take adventage of the following well-known triangle inequality
\begin{equation}\label{eq1.10}
\left|z_1+z_2\right|\leq\left|z_1\right|+\left|z_2\right|\hspace{2cm}(z_1,z_2\in\mathbb{C})
\end{equation}
and the following well-known geometric series sums
\begin{equation}\label{eq1.11}
\sum_{n\geq1}r^{n-1}=\frac{1}{1-r}
\end{equation}
and
\begin{equation}\label{eq1.12}
\sum_{n\geq1}nr^{n-1}=\frac{1}{(1-r)^2}
\end{equation}
for $r\in(0,1)$
in the proof of our results.

\section{Main Results}
\setcounter{equation}{0}
In this section, firstly, we present the following lemma which will be required in order to derive our main results.
\begin{lemma}\label{Lemma}
Let $i\in\{1,2,\dots,d\}, \alpha_i>-1$ and $2\lambda\mu>1,$ where $$\lambda=(d+1)^{d+1}\text{ and } \mu=\prod_{i=1}^{d}\left(\alpha_i+1\right).$$ Then, the normalized hyper-Bessel function $z\mapsto f_{\alpha_d}(z)$ satisfies the next two inequalities:
\begin{equation}\label{eq2.1}
\left|f_{\alpha_d}(z)\right|\leq\frac{2\lambda\mu+1}{2\lambda\mu-1}
\end{equation}
and
\begin{equation}\label{eq2.2}
\left|f_{\alpha_d}^{\prime}(z)\right|\leq\frac{4\lambda^2\mu\left(\mu+1\right)-1}{\left(2\lambda\mu-1\right)^2}.
\end{equation}
\end{lemma}

\begin{proof}
By using the inequalities \eqref{eq1.8}, \eqref{eq1.9} and \eqref{eq1.10} we can write that
\begin{align*}
\left|f_{\alpha_d}(z)\right|&=\left|z+\sum_{n\geq1}\frac{(-1)^n}{n!\left(d+1\right)^{n(d+1)}\prod_{i=1}^{d}(\alpha_i+1)_n}z^{n(d+1)+1}\right|\\&\leq1+\sum_{n\geq1}\frac{1}{2^{n-1}\left(d+1\right)^{n(d+1)}\prod_{i=1}^{d}(\alpha_i+1)^n}\\&=1+\frac{1}{\lambda\mu}\sum_{n\geq1}\left[\frac{1}{2\lambda\mu}\right]^{n-1}
\end{align*}
for $z\in\mathcal{U}.$ Here, using the geometric series sum which is given by \eqref{eq1.11} we deduce
$$\left|f_{\alpha_d}(z)\right|\leq\frac{2\lambda\mu+1}{2\lambda\mu-1}.$$
Similarly, in order to prove the inequality \eqref{eq2.2} we can use the inequalities which are given by \eqref{eq1.8}, \eqref{eq1.9}  and \eqref{eq1.10}. Namely,
\begin{align*}
\left|f_{\alpha_d}^{\prime}(z)\right|&=\left|1+\sum_{n\geq1}\frac{(nd+n+1)(-1)^n}{n!\left(d+1\right)^{n(d+1)}\prod_{i=1}^{d}(\alpha_i+1)_n}z^{n(d+1)}\right|\\&\leq1+\sum_{n\geq1}\frac{nd+n+1}{2^{n-1}\left(d+1\right)^{n(d+1)}\prod_{i=1}^{d}(\alpha_i+1)^n}\\&=1+\frac{1}{\mu}\sum_{n\geq1}\frac{n}{(2\lambda\mu)^{n-1}}+\frac{1}{\lambda\mu}\sum_{n\geq1}\left(\frac{1}{2\lambda\mu}\right)^{n-1}.
\end{align*}
Now, if we consider the geometric series sums which are given by \eqref{eq1.11} and \eqref{eq1.12}, then we have
$$\left|f_{\alpha_d}^{\prime}(z)\right|\leq\frac{4\lambda^2\mu\left(\mu+1\right)-1}{\left(2\lambda\mu-1\right)^2}.$$
So, the proof is completed.
\end{proof}

Let $w(z)$ denote an analytic function in $\mathcal{U}$. It is important to mention here that the following well-known result is very useful for our main results:
$$\Re\Bigg\{\frac{1+w(z)}{1-w(z)}\Bigg\}>0, \text{ if and only if }\left|w(z)\right|<1, z\in\mathcal{U}.$$

Now, we give our first main result related to the quotient of normalized hyper-Bessel function and its partial sum.
\begin{theorem}\label{the.1}
	Let $n\in\mathbb{N}=\{1,2,\dots\}, i\in\{1,2,\dots,d\}, \alpha_i>-1$, the function $f_{\alpha_d}:\mathcal{U}\to\mathbb{C}$ be defined by \eqref{eq1.7} and its sequence of partial sum defined by
	\begin{equation}\label{eq2.3}
	\left(f_{\alpha_d}\right)_m(z)=z+\sum_{n=1}^{m}A_nz^{n(d+1)+1}.
	\end{equation}  If the inequqlity $\lambda\mu>\frac{3}{2}$ is valid, then the following two inequalities are valid for $z\in\mathcal{U}:$
\begin{equation}\label{eq2.4}
\real\left(\frac{f_{\alpha_d}(z)}{\left(f_{\alpha_d}\right)_m(z)}\right)\geq\frac{2\lambda\mu-3}{2}
\end{equation}
and
\begin{equation}\label{eq2.5}
\real\left(\frac{\left(f_{\alpha_d}\right)_m(z)}{f_{\alpha_d}(z)}\right)\geq\frac{2\lambda\mu-1}{2}
\end{equation}
\end{theorem}
\begin{proof}
	From the inequality \eqref{eq2.1} in Lemma \eqref{Lemma} we can write that
	\begin{equation}\label{eq2.6}
	\left|f_{\alpha_d}(z)\right|=\left|z+\sum_{n\geq1}A_nz^{n(d+1)+1}\right|\leq1+\sum_{n\geq1}\left|A_n\right|\leq\frac{2\lambda\mu+1}{2\lambda\mu-1}.
	\end{equation}
	The inequality \eqref{eq2.6} is equivalent to
		\begin{equation}\label{eq2.7}
		\frac{2\lambda\mu-1}{2}\sum_{n\geq1}\left|A_n\right|\leq1.
		\end{equation}
		In order to prove the inequality \eqref{eq2.4}, we consider the function $w(z)$ defined by
		$$\frac{1+w(z)}{1-w(z)}=\frac{2\lambda\mu-1}{2}\Bigg\{\frac{f_{\alpha_d}(z)}{\left(f_{\alpha_d}\right)_m(z)}-\frac{2\lambda\mu-3}{2}\Bigg\}.$$ The last equality is equivalent to
\begin{equation}\label{eq2.8}
\frac{1+w(z)}{1-w(z)}=\frac{1+\sum_{n=1}^{m}A_nz^{n(d+1)}+\frac{2\lambda\mu-1}{2}\sum_{n=m+1}^{\infty}A_nz^{n(d+1)}}{1+\sum_{n=1}^{m}A_nz^{n(d+1)}}.
\end{equation}
Therefore, we obtain 
$$w(z)=\frac{\frac{2\lambda\mu-1}{2}\sum_{n=m+1}^{\infty}A_nz^{n(d+1)}}{2+2\sum_{n=1}^{m}A_nz^{n(d+1)}+\frac{2\lambda\mu-1}{2}\sum_{n=m+1}^{\infty}A_nz^{n(d+1)}}$$
and
$$w(z)\leq\frac{\frac{2\lambda\mu-1}{2}\sum_{n=m+1}^{\infty}\left|A_n\right|}{2-2\sum_{n=1}^{m}\left|A_n\right|-\frac{2\lambda\mu-1}{2}\sum_{n=m+1}^{\infty}\left|A_n\right|}.$$
The inequality
\begin{equation}\label{eq2.9}
\sum_{n=1}^{m}\left|A_n\right|+\frac{2\lambda\mu-1}{2}\sum_{n=m+1}^{\infty}\left|A_n\right|\leq1
\end{equation}
implies that $\left|w(z)\right|\leq1.$ It suffices to show that the left hand side of \eqref{eq2.9} is bounded above by
$$\frac{2\lambda\mu-1}{2}\sum_{n=1}^{\infty}\left|A_n\right|,$$
which is equivalent to
$$\frac{2\lambda\mu-3}{2}\sum_{n=1}^{m}\left|A_n\right|\geq0.$$ The last inequality holds true under the condition $\lambda\mu>\frac{3}{2}$.

The proof of the result \eqref{eq2.5} would run parallel to those of the result \eqref{eq2.4}. In order to do this, consider the function $p(z)$ defined by
\begin{align*}
\frac{1+p(z)}{1-p(z)}&=\left(1+\frac{2\lambda\mu-1}{2}\right)\Bigg\{\frac{f_{\alpha_d}(z)}{\left(f_{\alpha_d}\right)_m(z)}-\frac{2\lambda\mu-1}{2}\Bigg\}\\&=\frac{1+\sum_{n=1}^{m}A_nz^{n(d+1)}-\frac{2\lambda\mu-1}{2}\sum_{n=m+1}^{\infty}A_nz^{n(d+1)}}{1+\sum_{n=1}^{\infty}A_nz^{n(d+1)}}.
\end{align*}
So, we get that
$$p(z)=\frac{-\frac{2\lambda\mu+1}{2}\sum_{n=m+1}^{\infty}A_nz^{n(d+1)}}{2+2\sum_{n=1}^{m}A_nz^{n(d+1)}-\frac{2\lambda\mu-3}{2}\sum_{n=m+1}^{\infty}A_nz^{n(d+1)}}$$
and
$$p(z)\leq\frac{\frac{2\lambda\mu+1}{2}\sum_{n=m+1}^{\infty}\left|A_n\right|}{2-2\sum_{n=1}^{m}\left|A_n\right|-\frac{2\lambda\mu-3}{2}\sum_{n=m+1}^{\infty}\left|A_n\right|}.$$
The inequality \eqref{eq2.9} implies that $\left|p(z)\right|\leq1.$ Since the left hand side of the inequality \eqref{eq2.9} is bounded above
$$\frac{2\lambda\mu-1}{2}\sum_{n=1}^{\infty}\left|A_n\right|,$$
the proof is completed.
\end{proof}

Our second main result is the following:
\begin{theorem}\label{the.2}
	Let $n\in\mathbb{N}=\{1,2,\dots\}, i\in\{1,2,\dots,d\}, \alpha_i>-1$, the function $f_{\alpha_d}:\mathcal{U}\to\mathbb{C}$ be defined by \eqref{eq1.7} and its sequence of partial sum defined by \eqref{eq2.3}.  If the inequality $\frac{4\lambda^2\mu^2-4\lambda^2\mu-8\lambda\mu+3}{4\lambda^2\mu+4\lambda\mu-2}>0$ is valid, then the next two inequalities hold true for $z\in\mathcal{U}:$  
	\begin{equation}\label{eq2.10}
	\real\left(\frac{f_{\alpha_d}^{\prime}(z)}{\left(\left(f_{\alpha_d}\right)_m(z)\right)^{\prime}}\right)\geq\frac{4\lambda^2\mu^2-4\lambda^2\mu-8\lambda\mu+3}{4\lambda^2\mu+4\lambda\mu-2}
	\end{equation}
	and
	\begin{equation}\label{eq2.11}
	\real\left(\frac{\left(\left(f_{\alpha_d}\right)_m(z)\right)^{\prime}}{f_{\alpha_d}^{\prime}(z)}\right)\geq\frac{4\lambda^2\mu^2-4\lambda\mu+1}{4\lambda^2\mu+4\lambda\mu-2}.
	\end{equation}
\end{theorem}
\begin{proof}
	From the inequality \eqref{eq2.2} in Lemma \eqref{Lemma} we can write that
	\begin{equation}\label{eq2.12}
	1+\sum_{n\geq1}(nd+n+1)\left|A_n\right|\leq\frac{4\lambda^2\mu(\mu+1)-1}{(2\lambda\mu-1)^2}.
	\end{equation}
	The inequality \eqref{eq2.12} is equivalent to
	\begin{equation}\label{eq2.13}
	\frac{(2\lambda\mu-1)^2}{4\lambda^2\mu(\mu+1)-1}\sum_{n\geq1}(nd+n+1)\left|A_n\right|\leq1.
	\end{equation}
	In order to prove the inequality \eqref{eq2.10}, we consider the function $h(z)$ defined by
	\begin{align*}
	\frac{1+h(z)}{1-h(z)}&=\frac{(2\lambda\mu-1)^2}{4\lambda^2\mu(\mu+1)-1}\Bigg\{\frac{f_{\alpha_d}^{\prime}(z)}{\left(\left(f_{\alpha_d}\right)_m(z)\right)^{\prime}}-\frac{4\lambda^2\mu^2-4\lambda^2\mu-8\lambda\mu+3}{4\lambda^2\mu+4\lambda\mu-2}\Bigg\}\\&=\frac{1+\sum_{n=1}^{m}(nd+n+1)A_nz^{n(d+1)}+\frac{(2\lambda\mu-1)^2}{4\lambda^2\mu+4\lambda\mu-2}\sum_{n=m+1}^{\infty}(nd+n+1)A_nz^{n(d+1)}}{1+\sum_{n=1}^{m}(nd+n+1)A_nz^{n(d+1)}}.
	\end{align*}
	As a result, we get 
	$$h(z)=\frac{\frac{(2\lambda\mu-1)^2}{4\lambda^2\mu+4\lambda\mu-2}\sum_{n=m+1}^{\infty}(nd+n+1)A_nz^{n(d+1)}}{2+2\sum_{n=1}^{m}(nd+n+1)A_nz^{n(d+1)}+\frac{(2\lambda\mu-1)^2}{4\lambda^2\mu+4\lambda\mu-2}\sum_{n=m+1}^{\infty}A_nz^{n(d+1)}}$$
	and
	$$h(z)\leq\frac{\frac{(2\lambda\mu-1)^2}{4\lambda^2\mu+4\lambda\mu-2}\sum_{n=m+1}^{\infty}(nd+n+1)\left|A_n\right|}{2-2\sum_{n=1}^{m}(nd+n+1)\left|A_n\right|-\frac{(2\lambda\mu-1)^2}{4\lambda^2\mu+4\lambda\mu-2}\sum_{n=m+1}^{\infty}\left|A_n\right|}.$$
	The inequality
	\begin{equation}\label{eq2.14}
	\sum_{n=1}^{m}(nd+n+1)\left|A_n\right|+\frac{(2\lambda\mu-1)^2}{4\lambda^2\mu+4\lambda\mu-2}\sum_{n=m+1}^{\infty}(nd+n+1)\left|A_n\right|\leq1
	\end{equation}
	implies that $\left|h(z)\right|\leq1.$ It suffices to show that the left hand side of \eqref{eq2.14} is bounded above by
	$$\frac{(2\lambda\mu-1)^2}{4\lambda^2\mu+4\lambda\mu-2}\sum_{n=1}^{\infty}(nd+n+1)\left|A_n\right|,$$
	which is equivalent to	
	$$\frac{4\lambda^2\mu^2-4\lambda^2\mu-8\lambda\mu+3}{4\lambda^2\mu+4\lambda\mu-2}\sum_{n=1}^{m}(nd+n+1)\left|A_n\right|\geq0$$ such that the last inequality is valid under the our hypothesis.
	
	In order to prove the inequality\eqref{eq2.11}, consider the function $k(z)$ defined by
	\begin{align*}
	\frac{1+k(z)}{1-k(z)}&=\left(1+\frac{(2\lambda\mu-1)^2}{4\lambda^2\mu+4\lambda\mu-2}\right)\Bigg\{\frac{\left(\left(f_{\alpha_d}\right)_m(z)\right)^{\prime}}{f_{\alpha_d}^{\prime}(z)}-\frac{4\lambda^2\mu^2-4\lambda\mu+1}{4\lambda^2\mu+4\lambda\mu-2}\Bigg\}\\&=\frac{1+\sum_{n=1}^{m}(nd+n+1)A_nz^{n(d+1)}-\frac{4\lambda^2\mu^2-4\lambda\mu+1}{4\lambda^2\mu+4\lambda\mu-2}\sum_{n=m+1}^{\infty}(nd+n+1)A_nz^{n(d+1)}}{1+\sum_{n=1}^{\infty}(nd+n+1)A_nz^{n(d+1)}}.
	\end{align*}
	Consequently, we have that 
	$$k(z)=\frac{-\frac{4\lambda^2\mu^2-4\lambda\mu+1}{4\lambda^2\mu+4\lambda\mu-2}\sum_{n=m+1}^{\infty}(nd+n+1)A_nz^{n(d+1)}}{2+2\sum_{n=1}^{m}(nd+n+1)A_nz^{n(d+1)}-\frac{4\lambda^2\mu^2-4\lambda\mu+1}{4\lambda^2\mu+4\lambda\mu-2}\sum_{n=m+1}^{\infty}A_nz^{n(d+1)}}$$
	and
	$$k(z)\leq\frac{\frac{4\lambda^2\mu^2-4\lambda\mu+1}{4\lambda^2\mu+4\lambda\mu-2}\sum_{n=m+1}^{\infty}(nd+n+1)\left|A_n\right|}{2-2\sum_{n=1}^{m}(nd+n+1)\left|A_n\right|-\frac{4\lambda^2\mu^2-4\lambda\mu+1}{4\lambda^2\mu+4\lambda\mu-2}\sum_{n=m+1}^{\infty}\left|A_n\right|}.$$
	The inequality
	\begin{equation}\label{eq2.15}
	\sum_{n=1}^{m}(nd+n+1)\left|A_n\right|+\frac{4\lambda^2\mu^2-4\lambda\mu+1}{4\lambda^2\mu+4\lambda\mu-2}\sum_{n=m+1}^{\infty}(nd+n+1)\left|A_n\right|\leq1
	\end{equation}
	implies that $\left|k(z)\right|\leq1.$ Since the left hand side of the inequality \eqref{eq2.15} is bounded above
	$$\frac{4\lambda^2\mu^2-4\lambda\mu+1}{4\lambda^2\mu+4\lambda\mu-2}\sum_{n=1}^{\infty}(nd+n+1)\left|A_n\right|,$$
	the proof is completed.
\end{proof}

\section{Applications}
\setcounter{equation}{0}
In this section we present some applications concerning our main results. As we mentioned in the section 1, there is a relationship between hyper-Bessel function and classical Bessel function $J_\nu$. Clearly, for $d=1$ and $\alpha_1=\nu$, it is known that the normalized hyper-Bessel function $J_{\alpha_d}(z)$ reduces to classical Bessel function of the first kind given by
\begin{equation*}
J_\nu(z)=\sum_{n\geq0}\frac{(-1)^n}{n!\Gamma\left(\nu+n+1\right)}\left(\frac{z}{2}\right)^{2n+\nu}.
\end{equation*}
Also, putting $d=1$ and $\alpha_1=\nu$ in \eqref{eq1.7} we have the normalized classical Bessel function $\phi_\nu(z)=2^{\nu}\Gamma(\nu+1)z^{1-\nu}J_\nu(z)$ which has the following infinite sum representation:
\begin{equation}\label{eq3.1}
\phi_\nu(z)=\sum_{n\geq0}\frac{(-1)^n}{4^{n}n!(\nu+1)_n}z^{2n+1}.
\end{equation}
 By using this relationship the following corollaries can be given. Setting $d=1$ and $\alpha_1=\nu$ in Theorem \eqref{the.1}and \eqref{the.2}, respectively, we get the following:

\begin{corollary}\label{cor3.1}
	Let the function
	$\phi_\nu:\mathcal{U}\to\mathbb{C}$ be defined by \eqref{eq3.1}. The following assertions are valid for $z\in \mathcal{U}$. If $\nu>-\frac{5}{8}$, then
	\begin{equation}\label{eq3.2}
	\real\left(\frac{\phi_\nu(z)}{(\phi_\nu)_m(z)}\right)\geq\frac{8\nu+5}{2}
	\end{equation}
	and
	\begin{equation}\label{eq3.3}
	\real\left(\frac{(\phi_\nu)_m(z)}{\phi_\nu(z)}\right)\geq\frac{8\nu+7}{2}.
	\end{equation}	
\end{corollary}

\begin{corollary}\label{cor3.2}
Let the function
$\phi_\nu:\mathcal{U}\to\mathbb{C}$ be defined by \eqref{eq3.1}. The following assertions are valid for $z\in \mathcal{U}$. If $\nu>\nu^*$, then 	
\begin{equation}\label{eq3.4}
\real\left(\frac{{\phi^{\prime}}_\nu(z)}{\left((\phi_\nu)_m(z)\right)^{\prime}}\right)\geq\frac{64\nu^2+32\nu-29}{80\nu+78}
\end{equation}
\begin{equation}\label{eq3.5}
\real\left(\frac{\left((\phi_\nu)_m(z)\right)^{\prime}}{{\phi^{\prime}}_\nu(z)}\right)\geq\frac{64\nu^2+112\nu+49}{80\nu+78},
\end{equation}
where $\nu^*\approx0.46807$.
 \end{corollary}

 It is well-known from \cite{Watson} that there are the following relationships between elementary trigonometric functions and classical Bessel function $J_\nu$ for some special values of $\nu$:
\begin{equation}\label{eq3.6}
 J_{\frac{1}{2}}(z)=\sqrt{\frac{2}{\pi {z}}}\sin{z}\text{ and }J_{\frac{3}{2}}(z)=\sqrt{\frac{2}{\pi {z}}}\left(\frac{\sin{z}}{z}-\cos{z}\right).
\end{equation}
As a result of the above relationships, one can easily obtained that
\begin{equation}\label{eq3.7}
 \phi_{\frac{1}{2}}(z)=\sin{z}\text{ and }\phi_{\frac{3}{2}}(z)=3\left(\frac{\sin{z}}{z^2}-\frac{\cos{z}}{z}\right).
\end{equation}
Also, taking $m=0$ in the partial sums of trigonometric functions given by \eqref{eq3.7} we have
\begin{equation}\label{eq3.8}
\left(\phi_{\frac{1}{2}}\right)_0(z)=\left(\phi_{\frac{3}{2}}\right)_0(z)=z.
\end{equation}

\begin{example}\label{exam3.1}
In view of the Corollary \eqref{cor3.1} we have\\
{\bf a.} If we take $\nu=\frac{1}{2}$ and $m=0$ in \eqref{eq3.2} and \eqref{eq3.3}, respectively, then 
$$\real\left(\frac{\sin z}{z}\right)\geq\frac{9}{2}$$
and
$$\real\left(\frac{z}{\sin z}\right)\geq\frac{11}{2}.$$
{\bf b.} If we take $\nu=\frac{3}{2}$ and $m=0$ in \eqref{eq3.2} and \eqref{eq3.3}, respectively, then 
$$\real\left(\frac{\sin z-z\cos z}{z^3}\right)\geq\frac{17}{6}$$
and
$$\real\left(\frac{z^3}{\sin z-z\cos z}\right)\geq\frac{57}{2}.$$
\end{example}

\begin{example}\label{exam3.2}
		In view of the Corollary \eqref{cor3.2} we have\\
		{\bf a.} If we take $\nu=\frac{1}{2}$ and $m=0$ in \eqref{eq3.4} and \eqref{eq3.5}, respectively, then 
		$$\real(\cos z)\geq\frac{3}{118}$$
		and
		$$\real\left(\frac{1}{\cos z}\right)\geq\frac{118}{3}.$$
		{\bf b.} If we take $\nu=\frac{3}{2}$ and $m=0$ in \eqref{eq3.4} and \eqref{eq3.5}, respectively, then 
		$$\real\left(\frac{2z^2\cos z+(z^3-2z)\sin z}{z^4}\right)\geq\frac{163}{198}$$
		and
		$$\real\left(\frac{z^4}{2z^2\cos z+(z^3-2z)\sin z}\right)\geq\frac{361}{198}.$$
\end{example}

\begin{remark}
If we consider $m=0$ in the inequality \eqref{eq2.10}, then we obtain $\real\left(f_{\alpha_d}^{\prime}(z)\right)>0.$ In view of  the famous Noshiro-Warchawski Theorem (see \cite{Goodman}) we have that  the normalized hyper-Bessel function $f_{\alpha_d}$ is univalent in $\mathcal{U}$ for $\frac{4\lambda^2\mu^2-4\lambda^2\mu-8\lambda\mu+3}{4\lambda^2\mu+4\lambda\mu-2}>0.$
\end{remark}

\begin{remark}
If we consider $m=0$ in the inequality \eqref{eq3.4}, then we obtain $\real\left(\phi_{\nu}^{\prime}(z)\right)>0.$ In view of  the famous Noshiro-Warchawski Theorem (see \cite{Goodman}) we have that  the normalized Bessel function $\phi_{\nu}$ is univalent in $\mathcal{U}$ for $\nu>\nu^*\approx0.46807$.
\end{remark}

\end{document}